\documentclass[a4paper,12pt]{article}
\usepackage[spanish,USenglish]{babel}
\usepackage{times}
\usepackage[spanish]{babel}
\usepackage[latin1]{inputenc}
\usepackage{fancyhdr}
\usepackage{marvosym}
\usepackage{hyperref}
\usepackage{graphicx}
\usepackage{amsmath,amssymb,latexsym,amsthm,amscd}
\usepackage{calligra}
\usepackage{color}
\usepackage[T1]{fontenc}
\usepackage{lmodern}
\usepackage[latin1]{inputenc}
\usepackage{fancyhdr}
\usepackage{anysize}
\marginsize{3.0 cm}{3.0cm}{1.0 cm}{2.5 cm}
\usepackage[hmarginratio=1:1,top=25mm,bottom=20mm,left=25mm,right=25mm,columnsep=20pt]{geometry}

\usepackage{amsmath}
\usepackage{amsfonts}
\usepackage{amssymb}
\usepackage{graphicx}
\newtheorem{teorema}{Teorema}
\newtheorem{definicion}{Definici\'on}
\newtheorem{proposicion}{Proposici\'on}
\def\R{\mathop{I\!\!R}\nolimits}
\def\N{\mathop{I\!\!N}\nolimits}

\newtheorem{com}{Comentario}

	\title{	\large \bf Estimación de la dimensión de Hausdorff por dimensión de conteo por cajas y dimensión de información \\
		\vspace{0.4cm}
	}

\author{José Luis Ponte Bejarano.\thanks{Universidad Tecnológica del Perú, Trujillo - Perú. ({\tt c24064@utp.edu.pe}).}\and
	Juan Carlos Ponte Bejarano.\thanks{Universidad Nacional de Trujillo, Trujillo - Perú. ({\tt jponte@unitru.edu.pe}).}
	\and
	Alexis Rodriguez Carranza.\thanks{Universidad Privada del Norte, Trujillo - Perú. ({\tt alexis.rodriguez@upn.edu.pe}).}\and
	Eddy Cristiam Mirando Ramos.\thanks{Universidad Nacional de Trujillo, Trujillo - Perú. ({\tt emirandar@unitru.edu.pe}).} 
}
\begin{document}
 \maketitle
 

\selectlanguage{spanish}
\begin{center}
	{\bf Resumen}
\end{center}
{\it \noindent 
En el presente artículo se estima la dimensión de Hausdorff mediante la dimensión de conteo por cajas y la dimensión de información. Se demuestra que la primera es una cota superior para la dimensión de Hausdorff, mientras que la segunda es una cota inferior para la dimensión de conteo por cajas y que la dimensión de Información es cota superior para la dimensión de Hausdorff. Además, se estima la dimensión del atractor de Henon por dimensión de conteo por cajas.}

{\small {\bf Palabras clave}. Dimensión de la caja, dimensión de Hausdorff, Dimensión de Información.} 
 
 \selectlanguage{English}
\begin{center}
	{\bf Abstract}
\end{center}
{\it \noindent In this article, the Hausdorff dimension is estimated using the box-counting dimension and the information dimension. It is shown that the former is an upper bound for the Hausdorff dimension, while the latter is a lower bound for the box-counting dimension, and that the information dimension is an upper bound for the Hausdorff dimension. Additionally, the dimension of the Henon attractor is estimated using the box-counting dimension.}
 
{\small {\bf Keywords}. Box dimension, Hausdorff dimension, Information dimension.}
\selectlanguage{spanish}
\section{Introducion}
Una noción intuitiva de dimensión de un conjunto $A$ en $\R^{n}$ es como un exponente que expresa la escala del \textit{volumen} con su \textit{tamaño}:
$$\mbox{Volumen}\approx\mbox{tamaño}^{\mbox{dimensión}}$$
Aquí el \textit{volumen} puede expresar el volumen mismo, masa o incluso una medida de contenido de información, y \textit{tamaño} es una distancia lineal. Por ejemplo, \textit{volumen} de un cuerpo en dimensión 2 sería dado por su área, y el \textit{tamaño} sería dado por su diámetro \cite{Alexis}. Así, la dimensión es expresado como:
\begin{equation}\label{intr}
\mbox{Dimensión}=\lim_{\mbox{tamaño}\rightarrow 0}\frac{\log(\mbox{volumen})}{\log(\mbox{tamaño})},
\end{equation}
donde el límite indica que la dimensión es una cantidad local, y que para definir una cantidad global es preciso algún tipo de promedio \cite{Theiler}. Esta noción de dimensión, que se basa en el límite de un parámetro de tamaño pequeño que tiende a cero, aparece en todos lo conceptos de dimensión de un conjunto. 

La dimensión de Hausdorff, que desde el punto de vista teórico es la más importante y que tiene la ventaja de estar definida para cualquier conjunto, se define en términos de la medida de Hausdorff (\cite{Hausdorff},\cite{Falconer}) que teóricamente no es dificil de manipular. Esta dimensión tiene muchas propiedades matemáticamente útiles, pero en la práctica no es fácil de calcular o estimar numéricamente.

La dimensión de la caja llamada también de conteo por cajas, que tiene una forma similar a la igualdad (\ref{intr}), es una de las dimensiones más utilizadas debido a que es viable numéricamente y relativamente sencilla de calcular (\cite{Theiler},\cite{Sauer}). Además de este tipo de dimensión, otra es la dimensión de información, que indica el número de bits de información por simbolo para especificar la posición de un punto en un conjunto. Para definir este tipo de dimensión se emplea un concepto proporcionado por Claude Shannon, la entropía de la información \cite{Shannon}. Tanto la dimensión de la caja como la dimensión de información pueden ser útiles para estimar la dimensión de Hausdorff.

Tomando presente lo mencionado en el párrafo anterior, en este trabajo se usa la dimensión de la caja y la dimensión de información para estimar la dimensión de Hausdorff. Para ello, se mencionan tres resultados: el primero indica que la dimensión de la caja es una cota superior para la dimensión de Hausdorff; y el segundo indica que la dimensión de información es cota inferior para la dimensión de la caja. Para la demostración del segundo resultado se presenta una versión propia. El tercer resultado, indica que la dimensión de información es cota superior para la dimensión de Haurdorff. Finalmente, para ver cómo funcionan la dimensión de la caja en el cálculo de la dimensión de un conjunto, esta se usa para calcular o estimar la dimensión del atractor de Henon.

\section{Dimensión de Hausdorff y de la Caja}
\subsection{Dimensión de Haurdorff}
\begin{definicion}
Sean $A$ y $U$ subconjuntos no vacíos de $\R^{n}$. El diámetro de $U$, denotado por $\vert U\vert$, es definido por $$\vert U\vert := \sup\{\Vert x-y\Vert :\quad x,y\in U \}$$
Si $\{U_{i}\}$ es una colección contable o finita de subconjuntos de $\R^{n}$ de diámetro menor o igual a $\delta$ que cubre a $A$, es decir $A\subset\bigcup^{\infty}_{i=1}U_{i}$ con $0\leq\vert U_{i}\vert\leq\delta$ para cada $i$, decimos que $\{ U_{i}\}$ es un $\delta-$cubrimiento de $A$.
\end{definicion}
\begin{definicion}\label{Hausdorff}
Sea $A$ un subconjunto de $\R^{n}$ y $s$ un número real no negativo. Para cada $\delta>0$ se define:
\begin{equation}\label{eq2}
H_{\delta}^{s}(A):=\inf\{\sum_{i=1}^{\infty}\vert U_{i}\vert^{s}:\quad\{ U_{i}\}\quad \textsl{es un $\delta -$cubrimiento de $A$} \}
\end{equation}
Además, definimos:
\begin{equation}
H^{s}(A):=\lim_{\delta\rightarrow 0^{+}}H_{\delta}^{s}(A).
\end{equation}
$H^{s}(A)$ es llamada la medida de Hausdorff de $A$ de dimensión $s$.
\end{definicion} 

\begin{com}\label{com1} Apartir de la definición \ref{Hausdorff} se puede observar lo siguiente:\\
\begin{itemize}
\item[a)]Si $t>s$ y $\{ U_{i}\}$ es cualquier $\delta-$cubrimiento de $A$, entonces:
$$\sum_{i=1}^{\infty}\vert U_{i} \vert^{t}\leq\delta^{t-s}\sum_{i=1}^{\infty}\vert U_{i}\vert^{s}. $$
Tomando ínfimos respecto a los $\delta-$cubrimientos de $A$, se tiene:
$$\inf\{\sum_{i=1}^{\infty}\vert U_{i}\vert^{t}\}\leq\delta^{t-s}\inf\{\sum_{i=1}^{\infty}\vert U_{i}\vert^{s}\}.$$
De lo cual, por la ecuación (\ref{eq2}) de la definición \ref{Hausdorff}, se tiene:
$$H_{\delta}^{t}(A)\leq\delta^{t-s}H_{\delta}^{s}(A).$$
Si $H^{s}(A)<\infty$ y $\delta\rightarrow 0$, se tiene: $H^{t}(A)=0$ para $t>s$. 
\item[b)]Si $t<s$ y $\{ U_{i}\}$ es cualquier $\delta-$cubrimientos de $A$, entonces:
$$\sum_{i=1}^{\infty}\vert U_{i} \vert^{s}\leq\delta^{s-t}\sum_{i=1}^{\infty}\vert U_{i}\vert^{t}. $$
Tomando ínfimos respecto a los $\delta-$cubrimientos de $A$, se tiene:
$$\inf\{\sum_{i=1}^{\infty}\vert U_{i}\vert^{s}\}\leq\delta^{s-t}\inf\{\sum_{i=1}^{\infty}\vert U_{i}\vert^{t}\}.$$
De lo cual, por la ecuación (\ref{eq2}) de la definición \ref{Hausdorff}, se tiene:
$$H_{\delta}^{s}(A)\leq\delta^{s-t}H_{\delta}^{t}(A),$$
de donde:
$$H_{\delta}^{s}(A)\delta^{t-s}\leq H_{\delta}^{t}(A).$$
Si $H^{s}(A)<\infty$ y $\delta\rightarrow 0$, se tiene: $H^{t}(A)=\infty$ para $t<s$.
\end{itemize}
\end{com}
Tomemando presente lo indicado en el comentario \ref{com1}, establecemos la siguiente definición.

\begin{definicion}\label{dimH}
Sea $A$ un subconjunto de $\R^n$. La dimensión de Hausdorff de $A$, denotada por $\dim_{H}(A)$, está definida por:
$$\dim_{H}(A):=\inf\{s\geq 0: H^{s}(A)=0 \}=\sup\{s\geq 0 :H^{s}(A)=\infty \} $$
Así:
\begin{eqnarray*}
  H^{s}(A)=\left\{\begin{array}{l}
         \displaystyle\infty\mbox{  };\mbox{ si }0\leq s<\dim_{H}(A),\\
         \\
           0\mbox{  };\mbox{ si }s>\dim_{H}(A).\\
          \end{array}
   \right.
\end{eqnarray*}
\end{definicion}
Luego, de la definición \ref{dimH}, podemos obtener una representación gráfica de $H^{s}(A)$ y $\dim_{H}(A)$ (ver Figura \ref{Fig1}):
\begin{figure}[h]
  \begin{center}
 \scalebox{0.3}{\includegraphics{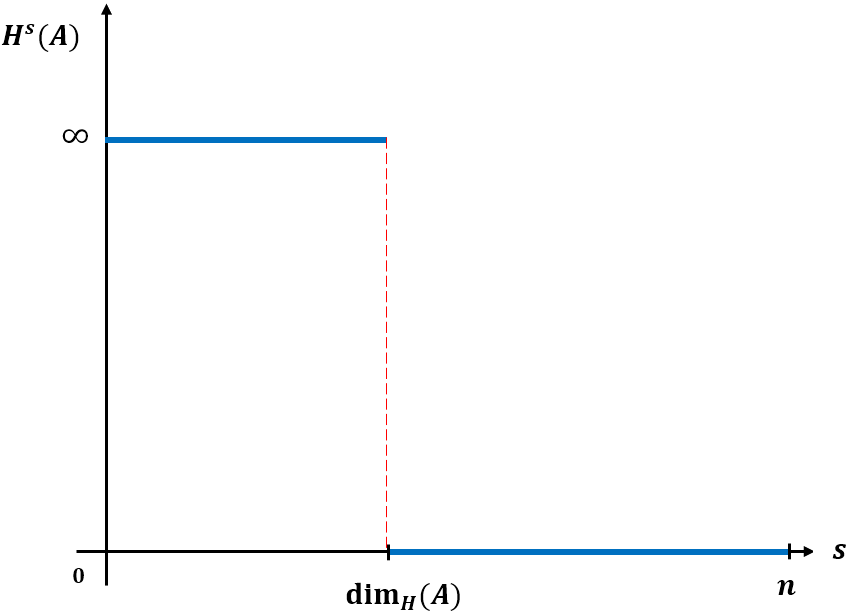}}
  \end{center}
 \caption{Grafica de $H^{s}(A)$ vs $s$.}
 \label{Fig1}
 \end{figure}

Algunas propiedades interesantes de la dimensión de Hausdorff desde el punto de vista teórico son la proposición \ref{propC1} y \ref{propC2}. La primera indica que la dimensión de Hausdorff de la imagen de un conjunto bajo una función de Lipschitz siempre es acotada por la dimensión de Hausdorff de ese conjunto. Mientras que la segunda indica que para determinar la dimensión de Hausdorff de un conjunto $A$, es suficiente conocer la mayor  dimensión de Haurdorff de algún conjunto que forma parte del cubrimiento del conjunto $A$. Estas proposiciones se indican a continuación.
\begin{proposicion}\label{propC1}
Sea $A$ un subconjunto de $\R^{n}$ y $f:\R^{n}\rightarrow\R^{n}$ una función de Lipschitz con constante de Lipschitz $c>0$. Entonces,
\begin{itemize}
\item[1)] Para cada $s>0$: $$\qquad H^{s}(f(A))\leq c^{s}H^{s}(A)$$
\item[2)] $\dim_{H}(f(A))\leq \dim_{H}(A)$
\end{itemize} 
\end{proposicion}
{\bf Prueba:}
\begin{itemize}
\item[1)] Sea $\{U_{i}\}$ un $\delta-$cubrimiento de $A$. Como $$|f(A\cap U_{i})|\leq c|A\cap U_{i}|\leq c|U_{i}|\leq c\delta,$$ entonces $\{f(A\cap U_{i})\}$ es $\varepsilon-$cubrimiento de $f(A)$, donde $\varepsilon =c\delta$. Así: $$\sum_{i=1}^{\infty}|f(A\cap U_{i})|^{s}\leq c^{s}\sum_{i=1}^{\infty}|U_{i}|^{s}.$$
Luego, $$H^{s}_{\varepsilon}(f(A))\leq c^{s}H^{s}_{\delta}(A).$$
Haciendo $\delta\rightarrow 0$, entonces $\varepsilon\rightarrow 0$ y  por tanto: $$H^{s}(f(A))\leq c^{s}H^{s}(A).$$
\item[2)] Si $s>\dim_{H}(A)$, entonces por la parte $1)$ se tiene que: $$H^{s}(f(A))\leq c^{s}H^{s}(A)=0.$$
De esto se tiene que: $\dim_{H}(f(A))\leq s$ para todo $s>\dim_{H}(A)$. Por tanto, $$\dim_{H}(f(A))\leq \dim_{H}(A).$$
\end{itemize}
\begin{proposicion}\label{propC2}
Sea $\{A_{j} \}_{j\in\N}$ una sucesión de conjuntos en $\R^{n}$ y $A=\displaystyle\bigcup_{j\in\N}A_{j}$. Entonces, $$\dim_{H}(A)=\sup_{j\in\N}\dim_{H}(A_{j}).$$
\end{proposicion}
{\bf Prueba:}\\
Como $A=\displaystyle\bigcup_{j\in\N}A_{j}$, entonces:
\begin{equation}\label{ec1Pro2}
\dim_{H}(A)\geq\sup_{j\in\N}\dim_{H}(A_{j})
\end{equation}
Por otro lado, si $s>\displaystyle\sup_{j\in\N}\dim_{H}(A_{j})$ entonces $$H^{s}(A)\leq \sum_{j\in\N}H^{s}(A_{j})=0,$$
de donde se tiene que $s\geq \dim_{H}(A)$. Como $s$ es arbitrario,
\begin{equation}\label{ec2Pro2}
\dim_{H}(A)\leq\sup_{j\in\N}\dim_{H}(A_{j})
\end{equation} 
De (\ref{ec1Pro2}) y (\ref{ec2Pro2}) se tiene: $$\qquad\dim_{H}(A)=\sup_{j\in\N}\dim_{H}(A_{j}).$$

A pesar de que la dimensión de Haurdorff tiene propiedades interesantes, como las proposiciones \ref{propC1} y \ref{propC2}, y otras que pueden ser revisadas en \cite{Falconer} y \cite{Hausdorff}, esta tiene desventajas. La principal desventaja, es la dificultad de su implementación numérica. Ello se debe al hecho de tener que calcular el ínfimo considerando todos los posibles cubrimientos. Una alternativa a ello, es relajar la búsqueda del ínfimo considerando todos los cubrimientos, para fijar una y refinarla. El valor que será obtenido es una cota superior para la dimensión de Hausdorff. Esta técnica es viable numéricamente y es conocida como la dimensión de la caja, la que se indica a continuación.

\subsection{Dimensión de la Caja}
\begin{figure}[h]
  \begin{center}
 \scalebox{0.25}{\includegraphics{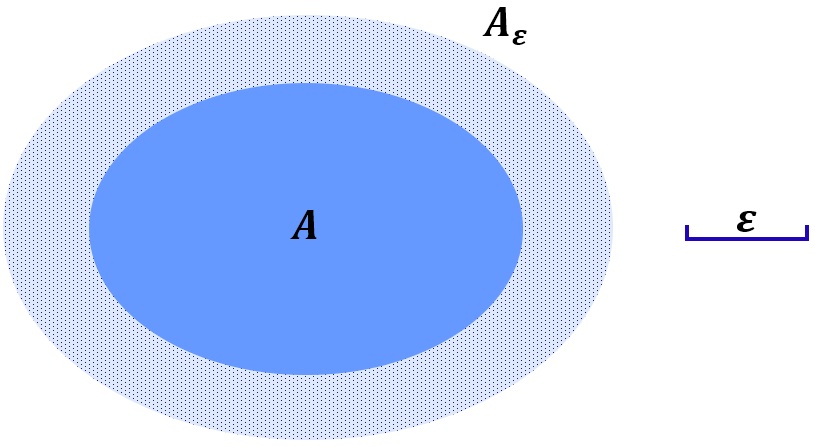}}
  \end{center}
 \caption{El conjunto $A$ y su $\varepsilon-$vecindad $A_{\varepsilon}$.}
 \label{Fractal-F1}
 \end{figure}

\begin{definicion}
Sea $A$ un subconjunto  de $\R^n$. Para cualquier $\varepsilon > 0$, se define la $\varepsilon-$vecindad o $\varepsilon-$cuerpo paralelo de $A$, denotado por $A_{\varepsilon}$, por:
$$A_{\varepsilon}:=\{x\in\R^{n}:\parallel x-a \parallel\leq\varepsilon\mbox{ para algún } a\in A \}.$$
Denotemos por $Vol\left( A_{\varepsilon}\right)$ el volumen $n-$dimensional de $A_{\varepsilon}$.
\end{definicion}
La $\varepsilon-$vecindad $A_{\varepsilon}$ representa el conjunto de puntos de $\R^{n}$ con distancia menor e igual a $\varepsilon$ de $A$. Una representación gráfica de este conjunto se muestra en la Figura \ref{Fractal-F1}.

\begin{figure}[h]
  \begin{center}
 \scalebox{0.43}{\includegraphics{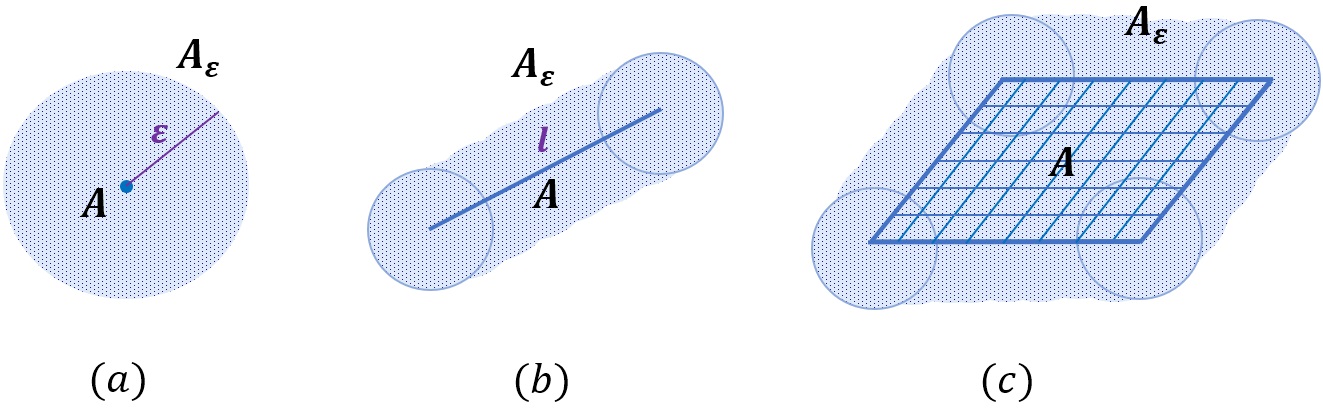}}
  \end{center}
 \caption{El conjunto $A$ cuando es un punto, un segmento de longitud $l$ y un plano de área $a$; y sus respectivos $\varepsilon-$vecindad $A_{\varepsilon}$.}
 \label{Fractal-F2}
 \end{figure}

Observemos a continuación algunos ejemplos del $Vol\left( A_{\varepsilon}\right)$ y de $A_{\varepsilon}$ de un conjunto $A$. Consideremos el espacio $\R^{3}$:
\begin{itemize}
\item [a)] Si $A$ es un sólo punto $x\in\R^{3}$, es decir, $A=\{x\}$ entonces $A_{\varepsilon}$ es una bola de centro $x$ y radio $\varepsilon>0$ con volumen $Vol(A_{\varepsilon})=\displaystyle\frac{4}{3}\pi\varepsilon^{3}$. Ver la imagen ($a$) de la Figura \ref{Fractal-F2}.
\item [b)]Si $A$ es un segmento de longitud $l$ entonces $A_{\varepsilon}$ representa un sólido de forma cilíndrica con volumen $Vol(A_{\varepsilon})\approx \pi l \varepsilon^{2}$. Ver la imagen ($b$) de la Figura \ref{Fractal-F2}.
\item [c)]Y si $A$ es un conjunto plano de área $a$ entonces $A_{\varepsilon}$ es un engrozamiento de $A$ con $Vol(A_{\varepsilon})\approx 2a \varepsilon$. Ver la imagen ($c$) de la Figura \ref{Fractal-F2}.
\end{itemize}
En los casos a), b) y c) se tiene que $Vol(A_{\varepsilon})\approx c\varepsilon^{3-s}$, donde el entero $s$ es la dimensión de $A$; y así el exponente de $\varepsilon$ es indicativo de la dimensión. El coeficiente $c$ de  $\varepsilon^{3-s}$ es una medida de longitud, área o volumen, según corresponda.

Esta idea se extiende a dimensiones fraccionarias. Si $A$ es un subconjunto de $\R^{n}$ y, para algún $s$, $Vol(A_{\varepsilon})/\varepsilon^{n-s}$ tiende a un número positivo $c$ cuando $\varepsilon\rightarrow 0$, entonces se tiene lo siguiente:
$$Vol(A_{\varepsilon})\approx c\varepsilon^{n-s}.$$
Tomando logaritmos,
$$\log(Vol(A_{\varepsilon}))\approx \log(c)+\log(\varepsilon^{n-s}),$$
de donde
$$\frac{\log(Vol(A_{\varepsilon}))}{\log(\varepsilon)}\approx \frac{\log(c)}{\log(\varepsilon)}+n-s$$
Luego,
$$\lim_{\varepsilon\rightarrow 0}\frac{\log(Vol(A_{\varepsilon}))}{\log(\varepsilon)} = n-s,$$
de donde:
$$s=n+\lim_{\varepsilon\rightarrow 0}\frac{\log(Vol(A_{\varepsilon}))}{-\log(\varepsilon)}.$$
Este hecho genera la siguiente definición.

\begin{definicion}\label{dimB1}
 Sea $A$ un subconjunto no vacio de $\R^n$.  La dimensión de la caja de $A$, llamada también dimensión de conteo por cajas, denotada por $\dim_{B}(A)$, es definida por:
$$\dim_{B}(A):=n+\lim_{\varepsilon\rightarrow 0}\frac{\log (Vol(A_{\varepsilon}))}{-\log (\varepsilon)},$$
si el límite existe. Si no existe, se define la dimensión de la caja usando los límites superior e inferior.
\end{definicion}

Hay una definición equivalente y más intuitiva para la dimensión de la caja. Para deducirla imaginemos el intervalo $[0,1]$ e intentemos calcular el número de intervalos de tamaño $\varepsilon >0$ necesarios para cubrirlo. Si fueran $n$ intervalos, tendríamos
$$n\varepsilon\approx 1,$$
de donde $n\approx\varepsilon^{-\dim ([0,1])}$. Para el caso $2-$dimensional pensemos en el cuadrado de lado $1$. Si $n$ cuadrados de lado $\varepsilon>0$ fueran necesarios para cubrirlo, tendríamos
$$n\varepsilon^{2}\approx 1,$$
de donde $n\approx\varepsilon^{-2}$. Esto nos lleva a pensar que si $\R^{n}$ es dividido en una malla de $\varepsilon-$cubos y $n(\varepsilon)$ es el número de cubos que tienen intersección con $A$, es decir, necesarios para cubrir $A$ entonces
$$n(\varepsilon)\approx\varepsilon^{-\dim(A)},$$
de donde se tiene:
 $$\dim(A)=\displaystyle\lim_{\varepsilon\rightarrow 0}\frac{\log n(\varepsilon)}{-\log(\varepsilon)}$$
 Este hecho permite establecer la siguiente definición.

\begin{definicion}\label{dimB2}
Sea $A$ un subconjunto no vacio de $\R^{n}$. La dimensión de la caja de $A$, llamada también dimensión de conteo por cajas, denotada por $\dim_B(A)$, es definida por:
 $$\dim_{B}(A):=\lim_{\varepsilon\rightarrow 0}\frac{\log n(\varepsilon)}{-\log(\varepsilon)},$$
si el límite existe; y donde $n(\varepsilon)$ representa el número de cubos de lado $\varepsilon > 0$ que tienen intersección con $A$, es decir, el número de cubos necesarios para cubrir $A$.
\end{definicion}
\begin{com}
En la definición \ref{dimB2}, aunque $n(\varepsilon)$ representa el número de cubos de lado $\varepsilon>0$ necesarios para cubrir el conjunto $A$, $n(\varepsilon)$ también puede representar:
\begin{itemize}
\item[1)]El menor número de bolas cerradas de radio $\varepsilon>0$ necesarios para cubrir $A$.
\item[2)]El menor número de conjuntos de diámetro a lo más $\varepsilon>0$ necesarios para cubrir $A$.
\item[3)] El mayor número de bolas disjuntas de radio $\varepsilon>0$ con centro en $A$.
\end{itemize}
\end{com}
\begin{proposicion}
Las definiciones \ref{dimB1} y \ref{dimB2} son equivalentes, es decir,
$$ n+\lim_{\varepsilon\rightarrow 0}\frac{\log (Vol(A_{\varepsilon}))}{-\log (\varepsilon)}=\lim_{\varepsilon\rightarrow 0}\frac{\log n(\varepsilon)}{-\log(\varepsilon)}$$
\end{proposicion}
{\bf Prueba:}\\
Sea $A$ un subconjunto de $\R^n$. Si $A$ se puede cubrir por $n(\varepsilon)$ bolas de radio $\varepsilon < 1$, entonces $A_{\varepsilon}$ puede ser cubierto por bolas concéntricos de radio $2\varepsilon$. Por tanto,
$$Vol(A_{\varepsilon})\leq n(\varepsilon)\cdot C_{n}\cdot(2\varepsilon)^{n},$$
donde $C_{n}$ es el volumen de la bola unitaria en $\R^{n}$. Tomando logaritmos:
$$\frac{\log(Vol(A_{\varepsilon}))}{-\log(\varepsilon)} \leq \frac{\log(2^{n}\cdot C_{n})+n\log(\varepsilon)+\log(n(\varepsilon))}{-\log(\varepsilon)},$$
de donde:
$$\lim_{\varepsilon\rightarrow 0}\frac{\log(Vol(A_{\varepsilon}))}{-\log(\varepsilon)} \leq -n + \lim_{\varepsilon\rightarrow 0}\frac{\log (n(\varepsilon))}{-\log(\varepsilon)}$$
Entonces,
\begin{equation}\label{pro3Ec1}
n + \lim_{\varepsilon\rightarrow 0}\frac{\log(Vol(A_{\varepsilon}))}{-\log(\varepsilon)} \leq \lim_{\varepsilon\rightarrow 0}\frac{\log (n(\varepsilon))}{-\log(\varepsilon)}
\end{equation}
Por otro lado, si existiera $n(\varepsilon)$ bolas disjuntas de lado $\varepsilon < 1$ con centros en $A$, entonces:
$$n(\varepsilon)\cdot C_{n}\cdot\varepsilon^{n}\leq Vol(A_{\varepsilon})$$
Tomando logaritmos:
$$\frac{\log(n(\varepsilon))+\log C_{n}+n\log(\varepsilon)}{-\log(\varepsilon)} \leq  \frac{\log(Vol(A_{\varepsilon}))}{-\log(\varepsilon)},$$
de donde
$$\lim_{\varepsilon\rightarrow 0}\frac{\log (n(\varepsilon))}{-\log(\varepsilon)}-n \leq \lim_{\varepsilon\rightarrow 0}\frac{\log(Vol(A_{\varepsilon}))}{-\log(\varepsilon)}$$
Entonces,
\begin{equation}\label{pro3Ec2}
\lim_{\varepsilon\rightarrow 0}\frac{\log (n(\varepsilon))}{-\log(\varepsilon)} \leq n+ \lim_{\varepsilon\rightarrow 0}\frac{\log(Vol(A_{\varepsilon}))}{-\log(\varepsilon)}
\end{equation}
De (\ref{pro3Ec1}) y (\ref{pro3Ec2}), se tiene:
$$ n+\lim_{\varepsilon\rightarrow 0}\frac{\log (Vol(A_{\varepsilon}))}{-\log (\varepsilon)}=\lim_{\varepsilon\rightarrow 0}\frac{\log n(\varepsilon)}{-\log(\varepsilon)},$$
lo que prueba que las definiciones \ref{dimB1} y \ref{dimB2} son equivalentes.

El siguiente teorema establece que la dimensión de la caja es cota superior para la dimensión de Hausdorff.
\begin{teorema}\label{Relacion1}
Sea $A\subset \R^{n}$. Entonces, $\dim_{H}(A)\leq \dim_{B}(A)$
\end{teorema}
\begin{proof}
Si $A$ se puede cubrir por $n(\varepsilon)$ conjuntos de diámetro $\varepsilon >0$, entonces por la definición \ref{Hausdorff}:
$$H_{\varepsilon}^{s}(A)\leq n(\varepsilon). \varepsilon^{s}$$
Si $1<H^{s}(A)=\displaystyle\lim_{\varepsilon\rightarrow 0}H_{\varepsilon}^{s}(A)$ entonces $\log n(\varepsilon)+ s \log(\varepsilon)>0$, para $\varepsilon >0$ suficientemente pequeño. Así,
$$s\leq \lim_{\varepsilon\rightarrow 0}\frac{\log n(\varepsilon)}{-\log(\varepsilon)},$$
de donde se tiene:
$$\dim_{H}(A)\leq \dim_{B}(A).$$
\end{proof}

Con el teorema \ref{Relacion1} se tiene una primera estimación para la dimensión de Haurdorff. Ahora vemos otra manera de pensar cuando se habla de la dimensión de un conjunto. Para ello introducimos la dimensión de información \cite{Theiler}.

\section{Dimensión de información}
Otra manera de pensar en la dimensión de un conjunto es en términos de cuántos números reales se necesitan para especificar la posición de un punto en ese conjunto (\cite{Theiler}, \cite{Alexis}). Por ejemplo, la posición de un punto en una recta es determinada por un número real, la posición en un plano con dos coordenadas cartesianas y la posición en el espacio tridimensional con tres coordenadas. Aquí, la dimensión es algo que cuenta el número de grados de libertad. Sin embargo, para conjuntos más complicados que líneas, planos y volúmenes, este concepto de dimensión debe extenderse.

Una manera de extender este concepto es determinar no cuántos números reales sino cuántos bits de información se necesitan para especificar la posición de un punto en un conjunto con una precisión determinada. Por ejemplo, en un segmento de línea de longitud unitaria, se necesitan $k$ bits para especificar la posición de un punto con una precisión de $r=2^{-k}$. Para un cuadrado unitario, se necesitan $2k$ bits para lograr la misma precisión ($k$ bits para cada lado). Y de manera similar, se necesitan $3$k bits para un cubo tridimensional. 
En general, $S(r)=-d\log_{2}(r)$ bits de información son necesarios para especificar la posición de un punto en un $d-$cubo  unitario  con una precisión  $r$. Estos ejemplos llevan a la definición natural de dimensión de información de un conjunto.
\begin{definicion}
Sea $A$ un conjunto no vacío de $R^{m}$. La dimensión de información del conjunto A, denotada por $\dim_I(A)$, es definida por:
$$\dim_I(A):=\lim_{r\rightarrow 0}\frac{-S(r)}{\log_{2}(r)}$$
donde $S(r)$ es la información (en bits) necesarios para especificar la posición de un punto sobre el conjunto $A$ con una precisión $r$.
\end{definicion}
El cálculo de $S(r)$ fue dada por Claude Shannon \cite{Shannon}. El mostró que si se tiene un conjunto de $n-$eventos, cuyas probabilidades de ocurrencia son $p_{1},p_{2},\ldots,p_{n}$, entonces la cantidad de información promedio para especificar una de ellas es dada por:
\begin{equation}\label{H}
H(p)=-\sum_{i=1}^{n}p_{i}\log_{2}(p_{i})
\end{equation}
A la función $H$, definida por la expresión (\ref{H}), se la llama \textit{la entropía de Shannon o de información}.

Para calcular la dimensión de información de un conjunto $A$ se debe considerar una partición del conjunto en $n-$ cajas $B_{i}$ de diámetro $r>0$. Luego, la probabilidad de ocurrencia de un punto del conjunto $A$ en la caja $B_{i}$ es dada por $$p_{i}=\displaystyle\frac{\mu(B_{i})}{\mu(A)},$$ donde $\mu(B_{i})=$Nro de puntos de $A$ en $B_{i}$ y $\mu(A)=$Nro de puntos de $A$. Usando la fórmula de Shannon (\ref{H}), la información necesaria para especificar un punto en el conjunto $A$ con una precisón $r>0$ es dada por:
$$S(r)=-\sum_{i=1}^{n}p_{i}\log_{2}(p_{i}).$$
Utilizando esta expresión, la dimensión de información de un conjunto $A$ se expresa por:
$$\dim_I(A)=\lim_{r\rightarrow 0}\frac{\displaystyle\sum_{i=1}^{n}p_{i}\log_{2}(p_{i})}{\log_{2}(r)}.$$

Con esta términología podemos enunciar el siguiente teorema, que indica que la dimensión de información es cota inferior para la dimensión de la caja.
\begin{teorema}\label{Relacion2}
Sea $A$ un subconjunto no vacio de $\R^m$ y $X$ una variable aleatoria discreta definida en $A$ con $n-$eventos, cuyas probabilidades de ocurrecia son $p_{1},p_{2},\ldots,p_{n}$. Entonces, 
$$\dim_I(A)\leq \dim_{B}(A)$$
\end{teorema}
\begin{proof}
Probemos inicialmente que $H(p)\leq \log_{2}(n)$.
En efecto:\\
Sea $u(x)=1/n$ la función de probabilidad uniforme sobre $A$, y $p(x)$ la función de probabilidad para $X$, cuyos valores son $p_{1}, p_{2},\dots, p_{n}$. Entonces:
\begin{eqnarray*}
-\sum_{x} p(x)\log_{2}\frac{p(x)}{u(x)}&=&\sum_{x} p(x)\log_{2}\frac{u(x)}{p(x)}\\
&\leq &\log_{2}\left(\sum_{x} p(x)\frac{u(x)}{p(x)}\right)\\
&=&\log_{2}\left(\sum_{x} u(x) \right)\\
&=&\log_{2}(1)\\
&=& 0,
\end{eqnarray*}
de donde, se obtiene 
\begin{equation}\label{H1}
\displaystyle\sum_{x} p(x)\log_{2}\frac{p(x)}{u(x)}\geq 0.
\end{equation}
Además, 
\begin{eqnarray*}
\sum_{x} p(x)\log_{2}\frac{p(x)}{u(x)}&=& \sum_{x}p(x)\log_{2}p(x)-\log_{2}p(x)\log_{2}u(x)\\
&=& \log_{2}(n)-H(p).
\end{eqnarray*}
Usando la igualdad anterior y la desigualdad (\ref{H1}) se tiene $H(p)\leq \log_{2}(n)$. 
De esta última desigualdad se puede afirmar que: $$S(r)\leq \log_{2}(n),$$
donde $n$ representa el número de bolas $B_{i}$ de radio $r>0$ (o cajas de diámetro $r>0$) que cubren $A$. Luego, para $r<1$:
$$\frac{S(r)}{-\log_{2}(r)}\leq \frac{\log_{2}(n)}{-\log_{2}(r)}.$$
Luego, $$\lim_{r\rightarrow 0}\frac{S(r)}{-\log_{2}(r)}\leq \lim_{r\rightarrow 0}\frac{\log_{2}(n)}{-\log_{2}(r)}.$$
Por tanto, $$\dim_I(A)\leq \dim_{B}(A)$$
Esto concluye la demostración.
\end{proof}

Hasta este punto, por los teoremas \ref{Relacion1} y \ref{Relacion2}, se ha establecido que la dimensión de la caja es cota superior de la dimensión de Hausdorff; y que la dimensión de información es cota inferior para la dimensión de la caja, así que surge una pregunta natural: La dimensión de información, ¿será una cota inferior ó superior para la dimensión de Haurdorff? En general no es posible establecer una relación directa entre ambas dimensiones. Sin embargo, bajo ciertas condiciones es posible establece una relación entre ellas: que la dimensión de Información es cota superior para la dimensión de Hausdorff, obteniendo de esta manera una estimación para la dimensión de Hausdorff. Esta se menciona a continuación.
\begin{teorema}\label{Relacion3}
Sea $A$ un subconjunto no vacio de $\R^m$ y $X$ una variable aleatoria discreta de $A$ con n$-$eventos cuyas probabilidades de ocurrencia son $p_1,p_2,\dots,p_n$. Si además $$S(\varepsilon)=\log_{2}(n(\varepsilon)),$$
donde $n(\varepsilon)$ es el número de bolas o cajas de diámetro $\varepsilon >0$ necesarios para cubrir $A$, entonces:
$$\dim_{H}(A)\leq \dim_{I}(A)$$
\end{teorema}
\begin{proof}
Sea $A$ un subconjunto  de $\R^{m}$ y $n(\varepsilon)$ el número de bolas o cajas de diámetros $\varepsilon >0$ necesarios para cubrir $A$, entonces por el Teorema \ref{Relacion1},
\begin{equation}\label{ec11}
0<\log_2(n(\varepsilon))+s\log_2(\varepsilon),
\end{equation}
siendo $\varepsilon >0$ suficientemente pequeño.
Además, por hipótesis del Teorema \ref{Relacion3},
\begin{equation}\label{ec12}
S(\varepsilon)=\log_2(n(\varepsilon))
\end{equation}
De la desigualdad (\ref{ec11}) y la igualdad (\ref{ec12}), se tiene:
$$0<S(\varepsilon)+s\log_2(\varepsilon),$$
de donde,
\begin{eqnarray*}
-s\log_2(\varepsilon)&<& S(\varepsilon)\\
\rightarrow s &<&\frac{S(\varepsilon)}{-\log_2(\varepsilon)}
\end{eqnarray*}
Cuando $\varepsilon\rightarrow 0$, se tiene:
$$s\leq \lim_{\varepsilon\rightarrow 0}\frac{-S(\varepsilon)}{\log_2(\varepsilon)}$$
es decir: $$\dim_H(A)\leq \dim_I(A)$$
\end{proof}

Con el Teorema \ref{Relacion1} y el Teorema \ref{Relacion3} se tienen dos maneras de estimar la dimensión de Hausdorff, y con el Teorema \ref{Relacion2} se puede afirmar que la dimensión de información dá una mejor estimación a la dimensión de Hausdorff, en comparación con la dimensión de conteo por cajas. Sin embargo, a pesar de ello la dimensión de conteo por cajas en más manejable en comparación con la dimensión de información, ya que para calcular la dimensión de información de un conjunto es necesario conocer las probabilidad de ocurrencia de los subconjuntos que forman parte de la partición, lo cual no es una tarea sencilla (ver el Apendice I en \cite{Farmer}). En cambio, para el calculo de la dimensión de conteo por cajas se necesita generar la partición del conjunto y determinar el número de cajas que tienen intersección con el conjunto. Como muestra de ello, a continuación pasamos a estimar la dimensión del atractor de Henon usando la dimensión de conteo por cajas.

\section{Dimensión del atractor de Henon}

El sistema de Henon \cite{Henon}, que es un modelo que representa la turbulencia de fluidos desde el punto de vista de un sistema dinámico, está dado por el siguiente sistema discreto:
\begin{eqnarray}
   \left\{\begin{array}{l}
         x(n+1)=1-ax(n)^2+y(n),\\\\         
         y(n+1)=bx(n),\\
          \end{array}
   \right.
\end{eqnarray}

\begin{figure}[h]
  \begin{center}
  \scalebox{0.52}{\includegraphics{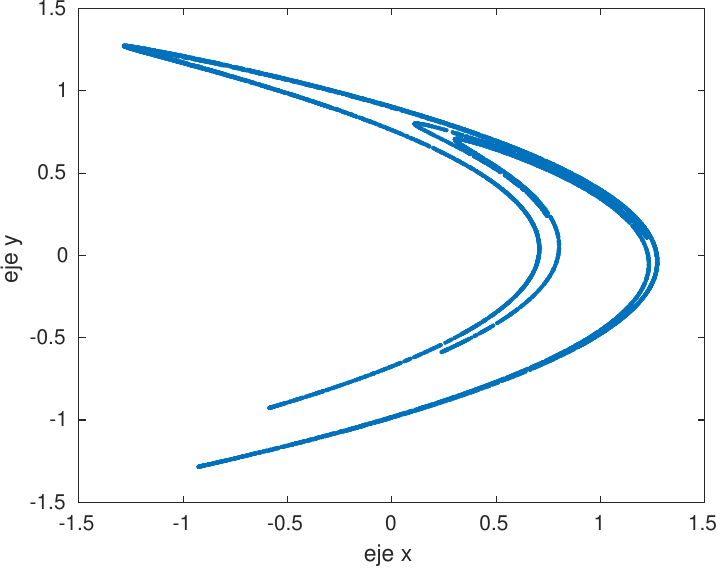}}
   \end{center}
  \caption{Atractor de Henon.}
  \label{henon-lorenz}
 \end{figure}

Cuando el valor de los parámetros es  $a=1$.4 y $b=0$.3 se genera un conjunto conocido como el atractor de Henon.
En la figura \ref{henon-lorenz} se muestra la representación geométrica de tal conjunto atractor.
 \\
Para estimar la dimensión del atractor de Henon se usará la dimensión de la caja. Para ello se debe contar el número de cajas o cuadros de lado $\varepsilon>0$ que lo intersectan.
\begin{figure}[h]
  \begin{center}
 \scalebox{0.27}{\includegraphics{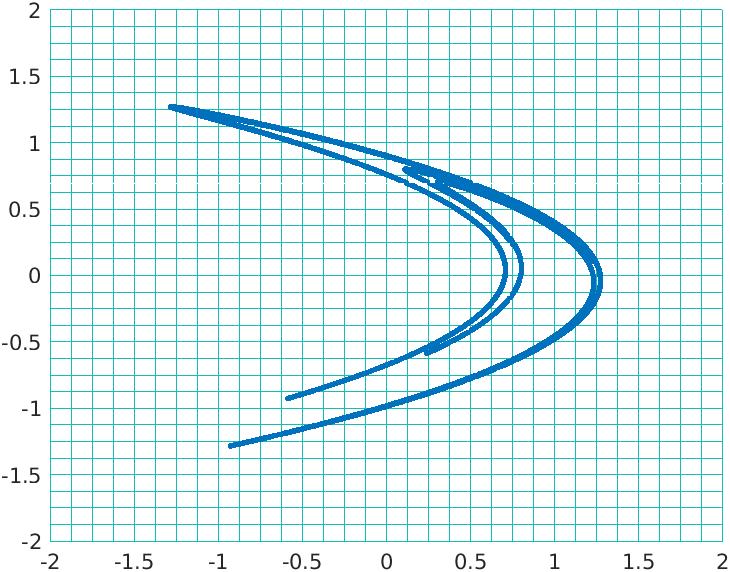}}
 \scalebox{0.27}{\includegraphics{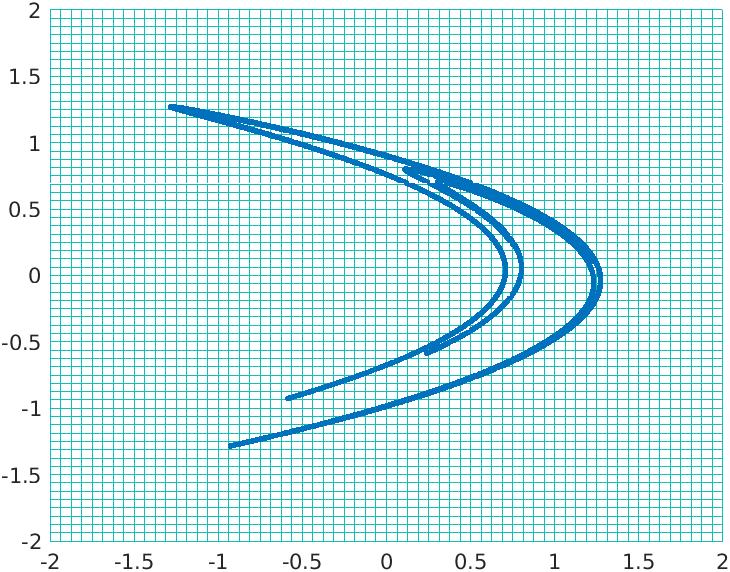}}
 \end{center}
 \caption{Cuadrículas con lados $\varepsilon=2^{-3}$ y $\varepsilon=2^{-4}$ respectivamente.}
 \label{dimB_henon}
 \end{figure}

En la figura \ref{dimB_henon} se observan dos cuadrículas, una con cuadros de lado $\varepsilon=2^{-3}$, y otra con $\varepsilon=2^{-4}$. Haciendo un conteo cuidadoso del número de cuadros que intersectan al atractor se tiene que $n(\varepsilon)=177$ para  $\varepsilon=2^{-3}$; y $n(\varepsilon)=433$ para  $\varepsilon=2^{-4}$. Haciendo el mismo procedimiento, se puede obtener que $n(\varepsilon)=1037$ para  $\varepsilon=2^{-5}$; $n(\varepsilon)=2467$ para  $\varepsilon=2^{-6}$ y  $n(\varepsilon)=5763$ para  $\varepsilon=2^{-7}$. Organizando esta información en un cuadro (cuadro \ref{cuadro1}) se tiene:
\begin{table}[htbp]
  \centering
     \begin{tabular}{|c|c|c|}
    \hline
    Cuadros de lados $\varepsilon>0 $ & Número de cuadros $n(\varepsilon)$\hbox\\
    \hline
    $2^{-3}$ & $177$ \hbox\\
    \hline
    $2^{-4}$ & $433$ \hbox\\
    \hline
    $2^{-5}$ & $1037$ \hbox\\
    \hline
    $2^{-6}$ & $2467$ \hbox\\
    \hline
    $2^{-7}$ & $5763$ \hbox\\
    \hline
    \end{tabular}%
     \caption{Número de cuadros $n(\varepsilon)$ de lado $\varepsilon $. }
     \label{cuadro1}
 \end{table}%
\\
Usaremos los datos del cuadro \ref{cuadro1} para estimar la dimensión del atractor de Henon. Como $\dim_{B}(A)=\displaystyle\lim_{\varepsilon\rightarrow 0}\frac{\log_{2} n(\varepsilon)}{\log_{2}(1/\varepsilon)}$, entonces podemos hacer una representación gráfica de $\log_{2}(n(\varepsilon))$ vs $\log_{2}(1/\varepsilon)$.
\begin{figure}[h]
  \begin{center}
 \scalebox{0.56}{\includegraphics{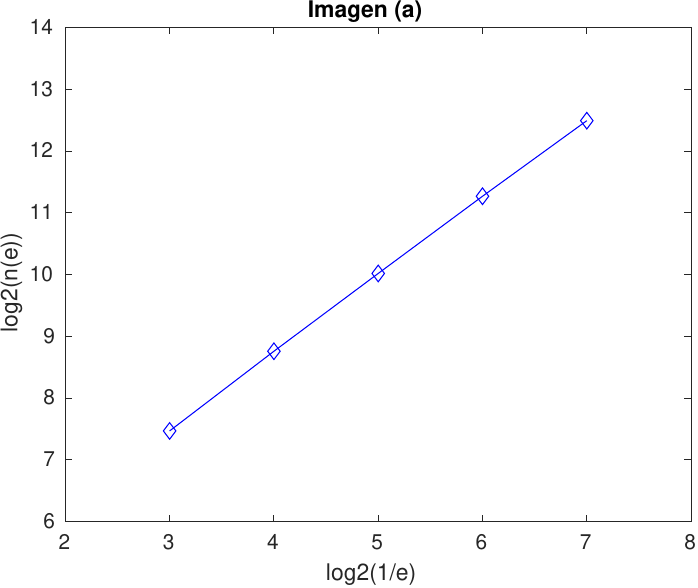}}
 \scalebox{0.56}{\includegraphics{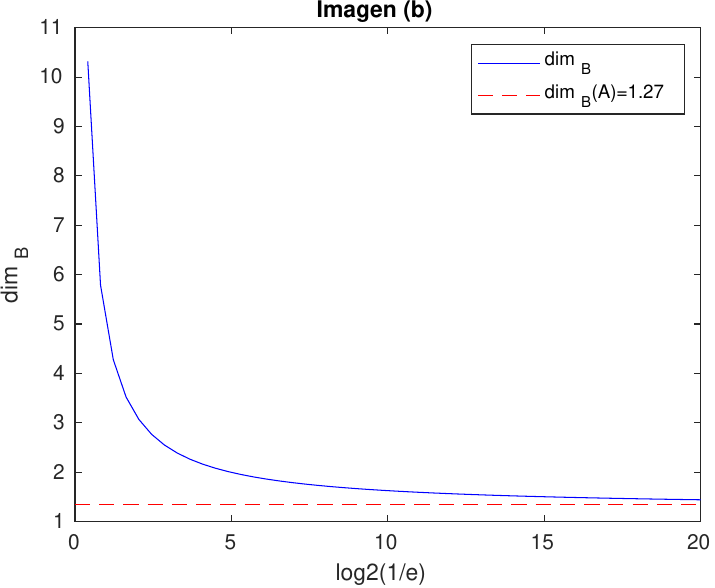}}
 \end{center}
 \caption{$\log_{2}(n(\varepsilon))$ vs $\log_{2}(1/\varepsilon)$ y $\dim_{B}$ vs $\log_{2}(1/\varepsilon)$.}
 \label{dimB_henon2}
 \end{figure}
\\
En la imagen (a) de la figura \ref{dimB_henon2}, se observa que la relación entre $\log_{2}(n(\varepsilon))$ y $\log_{2}(1/\varepsilon)$ sigue un comportamiento lineal. Entonces, para una cuadrícula de cuadros de lado $\varepsilon=2^{-b}$, con $b>0$ suficientemente grande
$$\log_{2}(n(\varepsilon))=\log_{2}\left(177\left(\frac{5763}{177}\right)^{\frac{b-3}{4}}\right)$$
Entonces, ver la imagen (b) de la figura \ref{dimB_henon2}, se tiene:
$$\dim_{B}(A)=\displaystyle\lim_{\varepsilon\rightarrow 0}\frac{\log_{2} n(\varepsilon)}{\log_{2}(1/\varepsilon)}=\lim_{b\rightarrow\infty}\frac{\log_{2}\left(177\left(\frac{5763}{177}\right)^{\frac{b-3}{4}}\right)}{\log_{2}(2^{b})}=1\mbox{.}27$$

Así, la dimensión del atractor de Henon, se estima por la dimensión de la caja en 1.27.

\section{Concluciones}
\begin{itemize}
\item[1)] La dimensión de la caja es cota superior para la dimensión de Hausdorff; mientras que la dimensión de información es cota inferior para la dimensión de conteo por cajas. Además, la dimensión de información es cota superior para la dimensión de Haurdorff.

\item[2)] De las dimensiones estudiadas, la dimensión de la caja es más manejable, ya que en el caso de la dimensión de Hausdoff al no conocer la forma de los delta cubrimientos se complica su utilización. Y en el caso de la dimensión de información, es necesario conocer las probabilidades de ocurrencia de los subconjuntos que forman parte de la partición, lo cual tampoco es una tarea sencilla.
\end{itemize}

\end{document}